\documentclass[11pt]{amsart}
\usepackage{amsfonts,latexsym,rawfonts,amsmath,amssymb}
\numberwithin{equation}{section}
\newcommand{\bthm}[2]{\vskip 8pt\bf #1\hskip 2pt\bf#2\it \hskip 8pt}
\newcommand{\ethm}{\vskip 8pt\rm}
\def\dint{\displaystyle{\int}}
\newcommand{\s}{\,\,\,\,}
\newcommand{\x}{\,\,\,\,}
\def\lan{\langle}
\def\ran{\rangle}

\begin{document}

\title[Moser-Trudinger inequalities for vector bundles]
{Moser-Trudinger inequalities of vector bundle over a compact
Riemannian manifold of dimension 2}
\author{Yuxiang Li, Pan Liu, Yunyan Yang}
\date{}
\maketitle
\begin{abstract}

Let $(M,g)$ be a 2-dimensional compact Riemannian manifold. In
this paper, we use the method of blowing up analysis to prove
several Moser-Trdinger type inequalities for vector bundle over
$(M,g)$. We also derive an upper bound of such inequalities under
the assumption that blowing up occur.

\end{abstract}

\section{\bf Introduction and Main results}

Let $(M, g)$ be a 2-dimensional compact Riemannian manifold. One
of Fontana's results (see \cite{Fo}) says
$$\sup_{\int_M|\nabla u|^2dV_g=1,\int_MudV_g=0}\dint_Me^{\alpha u^2}=\left\{
   \begin{array}{ll}
   <+\infty& if\s\alpha\leq 4\pi\\
   =+\infty& if\s\alpha>4\pi
   \end{array}\right.,
$$
which extends Trudinger and Moser's inequalities (see \cite{T},
\cite{M}). A weak form of the above inequality is
\begin{equation}\label{1.11}
\log \int_M e^u dV_g \leq \frac{1}{16\pi} \int_M | \nabla u
|^2dV_g +
 \int_M u dV_g + C
\end{equation}
for all $u \in H^{1,2}(M)$, where $C$ depends only on the geometry
of $M$ (see \cite{M}, \cite{Au}). The inequality (\ref{1.11})
 has been extensively used in many mathematical and
physical problems, for instance in the problem of prescribing
Gaussian curvature (\cite{Ch}, \cite{CY}, \cite{DJLW}), the mean
field equation and the abelian Chern-Simons model (\cite{DJLW3},
\cite{DJLW4}, \cite{J}), ect.

In this note we want to derive some new Moser-Trudinger type
inequalities. We will consider a smooth vector bundle $E$ with
metric $h$ and connection $\nabla$ over $M$. Throughout this
paper, we do not assume  $\nabla h=0$. To simplify the notations,
we write
\[
H_0 = \{ \sigma \in H^{1,2}(M, E): \nabla \sigma = 0 \},
\]
and $H_1 = H_0^{\perp}$, i.e.
\[
H_1 = \{ \sigma \in H^{1,2}(M, E): \int_M \lan \sigma, \zeta\ran
dV_g = 0, ~ \mbox{for all} ~ \zeta \in H_0 \}.
\]
We state our main results as follows:\\

\bthm{Theorem}{1.1} Let $(M,g)$ be a 2-dimensional compact
Riemannian manifold, $(E,h)$ be a smooth vector bundle over
$(M,g)$, $\nabla$ and $H_1$ be defined as above. Denote
\[
\mathcal{H}_1 = \{ \sigma \in H_1: \int_M \lan \nabla \sigma,
\nabla \sigma\ran d V_g = 1 \}.
\]
Then we have
\[
\sup_{\sigma \in \mathcal{H}_1} \int_M e^{ 4 \pi |\sigma |^2} d
V_g < \infty,
\]
where $|\sigma|^2=\langle\sigma,\sigma\rangle_h$ , and the
constant $4 \pi$ is sharp, which means that for any $ \alpha >
4\pi$,
\[
\sup_{\sigma \in \mathcal{H}_1} \int_M e^{ \alpha |\sigma |^2} d
V_g = \infty.
\] \ethm

As a corollary of Theorem 1.1, we have

\bthm{Corollary}{1.2} There exists a constant $C$ such that
\[
\int_M e^{|\sigma|} d V_g \leq C e^{\frac{1}{16\pi} \int_M \lan
\nabla \sigma, \nabla \sigma\ran dV_g}
\]
holds for all $\sigma \in H_1$.\ethm

We remark two special cases of Theorem 1.1: If $E$ is a trivial
bundle, ${e_i}$ is global basis of $E$, with
$$\lan e_i,e_j\ran_h=h_{ij}\x and \x \nabla e_i=0,$$
then we have the following:

\bthm{Corollary}{1.3} Let $(M,g)$ be a 2-dimensional compact
Riemannian manifold. Given an positive symmetric matrix $H$, we
denote
$$\mathcal{H}_1=\left\{U\in H^{1,2}(M,\mathbb{R}^n):
\int_M\nabla UH\nabla U^TdV_g=1, \int_M UdV_g=0\right\},$$
then we have
$$\sup_{U\in \mathcal{H}_1}
\int_Me^{4\pi UHU^T}dV_g<+\infty,$$ where $4\pi$ is sharp.\ethm

If $E$ is a trivial line bundle with $\lan e,e\ran_h=f(x)$,
Corollary 1.3  is exactly Yang's result\cite{Y}.\\

A complete analogue of Theorem 1.1 is the following:

\bthm{Theorem}{1.4.} Let $(M,g)$ be a 2-dimensional compact
Riemannian manifold, $(E,h)$ be a smooth vector bundle over
$(M,g)$, $\nabla$ and $H_1$ be defined as before. Denote
\[
\mathcal{H}_2 = \{ \sigma \in H_1: \int_M \lan \nabla \sigma, \nabla
\sigma\ran + |\sigma |^2) d V_g = 1 \}.
\]
Then we have
\[
\sup_{\sigma \in \mathcal{H}_2} \int_M e^{ 4 \pi |\sigma |^2} d
V_g < \infty,
\]
where the constant $4 \pi$ is sharp.\ethm

The proof of Theorem 1.1 is outlined as follows (the proof of
Theorem 1.4 is completely analogous to that of Theorem 1.1, so we
omit it). Let us define $J_{\alpha}(\sigma) = \int_M e^{\alpha
|\sigma|^2} dV_g$. We first show that the sup of $J_{\alpha}$ can
be attained in $\mathcal{H}_1$ if $\alpha < 4\pi$. So we can
choose $\alpha_k$ converging to $4\pi$ increasingly, and $\sigma_k
\in \mathcal{H}_1$ satisfying
\[
J_{\alpha_k}(\sigma_k) = \sup_{\sigma \in \mathcal{H}_1}
J_{\alpha_k}(\sigma).
\]
Denote $c_k = |\sigma|(x_k) = \max_{x \in M} |\sigma|(x)$. Passing
to a subsequence, we assume $ p =
\lim\limits_{k\rightarrow+\infty} x_k$. Without loss of
generality, we may assume blow-up occur, that is, $c_k \rightarrow
+\infty$. Take a local coordinate system $(\Omega, x)$ around $p$.
Using the idea in \cite{S}, we define a sequence of functions,
\[
\varphi_k (x) = 2 \alpha_k \lan\sigma_k(x_k), \sigma_k(x_k + r_k x) -
\sigma_k(x_k)\ran
\]
for some $r_k > 0$, where $x \in \Omega_k$ with
\[
\Omega_k = \{ x \in \mathbb{R}^2 ~ | ~ x_k + r_k x \in \Omega\}.
\]
We then prove that, for suitable $r_k$,
\[
\varphi_k \rightarrow - 2 \log ( 1 + \pi |x|^2)\quad  \mbox{in}
\quad C^2_{\mbox{loc}}(\mathbb{R}^2).
\]
Then we prove that
\[
c_k \sigma_k \rightarrow G \quad \mbox{in}\quad C^2(\Omega), \quad
\forall \Omega \subset \subset M \setminus \{p\},
\]
where $G$ is certain Green section. Finally, with the asymptotic
behavior of $\sigma_k$ described above, we can establish the
desired inequality and thus Theorem 1.1. In fact, we can give an
upper bound of the functional $J_{4\pi}$  in case that the
blow-up happens.\\

Though we mainly follow the ideas in \cite{L} and \cite{LL}, we
should point out that, in this paper, the convergence of
$c_k\sigma_k$ is derived differently from that of \cite{L} or
\cite{LL}: The key gradient in \cite{L} is the energy estimate
\begin{equation}{\label{1.1}}
\lim_{k\rightarrow+\infty}\dint_{\{u_k\leq\frac{c_k}{A}\}} |\nabla
u_k|^2dV_g=\frac{1}{A}.
\end{equation}
In this paper, though similar identity is used, the calculations
are based on the local Pohozaev identities. We thank Professor
Weiyue Ding who notice us possible application of the Pohozaev
identity when we study the extremal function for Fontana's
inequality on 4-dimensional manifold (see \cite{LY}). Moreover,
the method we get the upper bound of $J_{4\pi}$ is also new:
instead of capacity technique in \cite{L}, we use a result of
Carleson and Chang (\cite{CC}) as follows:

\bthm{Theorem}{A} {Let $B$  be unit ball in $\mathbb{R}^2$. Given
any sequence $u_k\in H^{1,2}_0(B)$, if
$u_k\rightharpoondown 0$, and $\int_{B_\delta}|\nabla u_k|^2dx\leq
1$, then we have
$$\limsup_{k\rightarrow +\infty}\int_{B}(e^{4\pi u_k^2}-1)dx\leq\pi e.$$}\ethm

The paper is organized as follows: In section 2, we settle some
notations for use later. In section 3 we prove that $4\pi$ is the
best constant. Section 4 is blowing up analysis. We will prove the convergence of
$c_k\sigma_k$ in section 5. Then we finish the
proof of our main theorem in section 6.

We hope our results can be a powerful tool in studying some
problems arising from geometry and mathematical physics. In a
forthcoming paper, we shall extend our results to high dimensional
case and find some geometrical and physical applications.

\section{\bf Preliminaries}
In this section, we clarify some notations.
 Take finite coordinate domains
$\{\Omega_k\}$ which cover M. Let $\sigma$ be a smooth section of
$E$, on each $\Omega_k$, we can set $\sigma=\sum_{i=1}^nu^ie_i$.
We define
$$||\sigma||_{H^{1,2}_k}^2=\int_M\sum_{i,j}(|\frac{\partial u^i}{\partial x^j}|^2+
\sum_k|u^i|^2)dx$$ and
$$||\sigma||_{H^{1,2}(M,E)}=\sum_{k}||\sigma||_{H^{1,2}_k}.$$
Clearly, such norm is equivalent to
$$||\cdot||=\left(\int_M(|\nabla \cdot|^2+|\cdot|^2)dV_g\right)^\frac{1}{2}.$$

Let $\sigma$ be a parallel section, i.e. $\sigma\in H_0$. Then we
have $||\sigma||_{H^{1,2}}\leq C||\sigma||_{L^2}$. By the
compactness of Sobolev embedding from $H^{1,2}$ into $L^2$, we
have

1. $H_0$ is a finite dimensional vector space. Throughout this
paper, we use $\zeta_1,\cdots,\zeta_m$ to denote an  orthogonal
basis of $H_0$.

2. Poicar\'{e} inequality holds on $H_1$, i.e. for any $\sigma\in
H_1$, we have
$$\int_M|\sigma|^2dV_g\leq C\int_M\lan\nabla\sigma,\nabla\sigma\ran dV_g,$$
then, on $H_1$, we can set
$$||\cdot||_{H^{1,2}_1}=\left(\int_M(|\nabla\cdot|^2)dV_g\right)^\frac{1}{2}$$
which is equivalent to the $H^{1,2}$ norm.\\

Throughout this paper,  we will use $(\Omega_p;x^1,x^2)$ to denote
 an isothermal coordinate system around $p$ with $p=(0,0)$. We can write the
 metric in such coordinate system as follows:
$$g=e^{f_p}(d(x^1)^2+d(x^2)^2)$$
with $f_p(0)=0$. Moreover, we always assume $e_1,e_2,\cdots,e_n$ to be an orthogonal
basis of $E$ in $\Omega_p$.

We should explain some notations involving $\nabla\Phi$. Locally,
when $\Phi=u^ke_k$ is a section of $E$, and $\nabla$ denotes the
connection of $E$, then
$$\lan\nabla\Phi,\nabla\Phi\ran= e^{-f_p}\sum_{i=1,2}
\lan\nabla_{\frac{\partial}{\partial x^i}}u^ke_k
,\nabla_{\frac{\partial}{\partial x^i}}u^ke_k\ran= \sum_{i=1,2}
\lan\nabla_{\frac{\partial}{\partial x^i}}u^ke_k
,\nabla_{\frac{\partial}{\partial x^i}}u^ke_k\ran_0.$$ When $\Phi$
is a function, then $\nabla\Phi$ is just the tangent vector
$\frac{\partial\Phi}{\partial x^1} \frac{\partial}{\partial x^1}+
\frac{\partial\Phi}{\partial x^2}\frac{\partial}{\partial x^2}$,
and
$$|\nabla_0\Phi|^2=|\frac{\partial\Phi}{\partial x^1}|^2+
|\frac{\partial\Phi}{\partial x^2}|^2,\quad |\nabla\Phi|^2=
e^{-f_p}|\nabla_0\Phi|^2.$$ When $\Phi=(u^1,\cdots,u^n)\in
H^{1,2}(\Omega)\times\cdots\times H^{1,2}(\Omega)$, we use
$|\nabla_0\Phi|^2$ to denote
$$\sum_{k=1}^n|\nabla_0 u^k|^2.$$

Throughout this paper we use $B_r$ to denote the following open
set of $\Omega_p$
$$\{(x^1,x^2):|x^1|^2+|x^2|^2<r^2\}.$$

Finally, we use $\Delta_0$ to denote the standard Laplacian on
$\mathbb{R}^2$
$$\frac{\partial^2}{\partial {x^1}^2}+\frac{\partial^2}{\partial{x^2}^2}.$$

\section{\bf The best constants}

The main task of this section is to prove the following:  $\sup_{{\sigma}\in \mathcal{H}_1}\int_M
e^{\alpha |\sigma|^2}dV_g< +\infty$ for any $\alpha<4\pi,$ and
$\sup_{{\sigma}\in \mathcal{H}_1}\int_{M}
  e^{\alpha |\sigma|^2}dV_g= +\infty$ for any $\alpha>4\pi.$
Firstly, we need a result in \cite{Ch} (cf. Theorem 2.50 in [Au]):

\bthm{Lemma}{3.1.} Let $B_r$ be
a ball in $\mathbb{R}^2$, then for any $\beta<2\pi$, we have
$$\sup_{\int_{B_r}(|\nabla u|^2+u^2)dx=1}\int_{B_r}e^{\beta u^2}dx<C(r).$$
\ethm

Then we have the following:

\bthm{Corollary}{3.2.} Let $B_r$ be a ball in $\mathbb{R}^2$, and
$u^1,\cdots,u^n\in H^{1,2}(B_r)$. Then we have
$$\sup_{\sum\limits_{i=1}^n\int_{B_r}(|\nabla u^i|^2+|u^i|^2)dx=1}\int_{B_r}e^{
\beta\sum\limits_{i=1}^n|u^i|^2}dx<C(r).$$
for any $\beta<2\pi$.\ethm

\proof We set $\lambda_i=\int_M(|\nabla u^i|^2+|u^i|^2)dx$, then
$\int_\Omega e^{\beta \frac{|u^i|^2}{\lambda_i}}dx<C$. Hence
$$\int_\Omega e^{\beta\sum\limits_{i=1}^n|u^i|^2}dx\leq\prod_{i=1}^n(\int_\Omega
e^{\beta \frac{|u^i|^2}{\lambda_i}}dx
)^{\lambda_i}\leq C^{\sum\limits_{i=1}^n\lambda_i}=C.$$

$\hfill\Box$\\

{\bf Remark 3.1.} Let $\Omega_p$ be a coordinate system around
$p$, and $B_r\subset\Omega_p$. Given a section $\sigma=u^ke_k$, we
set $U= (u^1,\cdots,u^n)$. Clearly, $|\nabla_0 U|^2+|U|^2\leq
C(\lan\nabla\sigma,\nabla\sigma\ran+|\sigma|^2)$. Hence, for any
$\alpha<\frac{2\pi}{C}$, we have
$$\sup_{\int_{B_r}(\lan\nabla\sigma,\nabla\sigma\ran+|\sigma|^2)
dV_g=1}\int_{B_r} e^{\alpha |\sigma|^2}dV_g<+\infty.$$
As an
consequence, we have  the following
  lemma:

\bthm{Lemma}{3.3.}  There exists a positive number $\alpha$ such
that $\sup_{{\sigma}\in \mathcal{H}_1}\int_M e^{\alpha
|\sigma|^2}dV_g<+\infty.$\ethm

Denote
$\tilde{\alpha}=\sup\{\alpha: \sup_{{\sigma}\in \mathcal{H}_1}\int_M e^{\alpha
  |\sigma|^2}dV_g<+\infty\}$. We shall prove that $\tilde{\alpha}=4\pi$.

\bthm{Lemma}{3.4.}  $\tilde{\alpha}\leq 4\pi$.\ethm

\proof Let $\Omega_p$ be  a coordinate domain. By Moser's result
\cite{M},  we can find a sequence $\{u_k\}\subset
H^{1,2}_0(B_\frac{1}{k})$ such that $\int_{B_\frac{1}{k}}|\nabla_0
u_k|^2dx=\int_{B_\frac{1}{k}}|\nabla u_k|^2dV_g=1$, and
$\int_{B_\frac{1}{k}}e^{(4\pi+\frac{1}{k})u_k^2}dV_g>k$ as
$k\rightarrow +\infty$. We set
$$\sigma_k=u_ke_1-\sum_{i=1}^m\lan u_ke_1,\zeta_i\ran_{L^2}\zeta_i\s\in H_1.$$
Since $||u_k||_{L^2}\rightarrow 0$, we get
$$\int_M\lan\nabla\sigma_k,\nabla\sigma_k\ran dV_g\rightarrow 1,\s {\rm and} \s
\int_Me^{(4\pi+\frac{1}{k})|\sigma_k|^2}dV_g\rightarrow\infty.$$

 $\hfill\Box$\\

Next, we prove an energy concentration phenomenon as follows:

\bthm{Lemma}{3.5.}  Let ${\sigma}_k\in \mathcal{H}_1$. If
\begin{equation}\label{3.1}\lim_{k\rightarrow+\infty}\int_Me^{\alpha
|\sigma_k|^2}dV_g=+\infty
\end{equation}
for all $\alpha>\tilde{\alpha}$, then passing to a
subsequence,${\sigma}_k\rightharpoondown 0$ in $H^{1,2}(M,E)$, and
$$\lan\nabla\sigma_k,\nabla\sigma_k\ran dV_g\rightharpoondown\delta_p$$
 for some $p\in M$.\ethm

{\bf Proof.} Without loss of generality, we assume
$$\begin{array}{lll}
 \sigma_k\rightharpoonup \sigma_0 &{\rm weakly\,\,\, in}& H^{1,2}(M,E)\\
 \sigma_k\rightarrow \sigma_0 &{\rm strongly\,\,\, in}& L^2(M,E).
\end{array}$$
We first claim that $\sigma_0=0$. Suppose not, we have
$$\int_M\lan\nabla(\sigma_k-\sigma_0),\nabla(\sigma_k-\sigma_0)\ran dV_g
 \rightarrow 1-\int_M\lan\nabla\sigma_0,\nabla\sigma_0\ran dV_g<1.$$
Hence we can find a $\alpha_1>\tilde{\alpha}$ such that
$\int_Me^{\alpha_1|\sigma_k|^2}dV_g$ is bounded, which
contradicts (3.1).

Secondly we show the concentration phenomenon. Suppose
$$\lan\nabla\sigma_k,\nabla\sigma_k\ran dV_g\rightharpoonup
 \mu\neq\delta_p,\quad\forall p\in M.$$
Taking $\eta\in C_0^\infty(B_\delta(p))$, $\eta\equiv 1$ on
$B_{\delta/2}(p)$, one can easily see that
$$\limsup_{k\rightarrow+\infty}\dint_{B_\delta(p)}\lan
\nabla(\eta\sigma_k), \nabla(\eta\sigma_k)\ran dV_g<1.$$ Hence
$e^{|\eta\sigma_k|^2}$ is bounded in $L^\alpha(M,E)$ for some
$\alpha>\tilde{\alpha}$. A covering argument implies that
$e^{|\sigma_k|^2}$ is bounded in $L^\alpha(M,E)$ for some
$\alpha>\tilde{\alpha}$, which contradicts (3.1).
$\hfill\Box$\\

\bthm{Corollary}{3.6.} $\tilde{\alpha}=4\pi$\ethm

{\bf Proof.} By the definition of $\tilde{\alpha}$, we can find a
sequence ${\sigma}_k\in \mathcal{H}_1$ such that
$$\lim_{k\rightarrow+\infty}\int_Me^{(\tilde{\alpha}+\frac{1}{k})|\sigma_k|^2}dV_g=+\infty.$$
Then, applying Lemma 3.5, we get
$$\lan\nabla\sigma_k,\nabla\sigma_k\ran dV_g\rightharpoondown\delta_p,\s
{\rm and} \s \dint_M|\sigma_k|^2dV_g\rightarrow 0.$$ Hence, for
any $\eta$ which is 0 near p, we have $\int(|\nabla(\eta
{\sigma}_k)|^2+|\eta\sigma_k|^2)dV_g \rightarrow 0$. Applying
Lemma 3.3, one can find a subsequence (still denoted by
$\sigma_k$) such that
$$\lim_{k\rightarrow+\infty}\int_{\Omega}e^{q |\eta\sigma_k|^2}dV_g=Vol(\Omega)$$
for any $\Omega\subset\subset M\setminus\{p\}$ and $q>0$.

In coordinate around $p$, we  set $\sigma_k=u_k^ie_i$. It is easy
to check that
$$\lim_{k\rightarrow+\infty}\dint_{B_{2r}}\sum_{i,j}
|\frac{\partial{(\eta'u_k^i)}}{\partial{x^j}}|^2dx=1,$$ where
$\eta'$ is a cut-off function which is 1 on $B_r$. Then,  similar
to the proof of Corollary 3.2, we can deduce from Moser's result
\cite{M} that
$$\lim_{k\rightarrow+\infty}\int_{B_r}e^{q |\sigma_k|^2}dx<+\infty$$
for any $q<4\pi$. Hence $\tilde{\alpha}\geq 4\pi$, which together
with Lemma 3.4 implies $\tilde{\alpha}=4\pi$.
$\hfill\Box$\\

In a similar way, we can prove the following

\bthm{Proposition}{3.7} For any $\alpha<4\pi$,
$$\sup_{\sigma\in \mathcal{H}_2}\int_Me^{\alpha |\sigma|^2}dV_g<+\infty,$$
and for any $\alpha>4\pi$,
$$\sup_{\sigma\in \mathcal{H}_2}\int_Me^{\alpha |\sigma|^2}dV_g=+\infty.$$\ethm

\section{\bf Blowing up analysis}

Let $\alpha_k$ be an increasing sequence which converges to
$4\pi$. In this section, we shall consider a sequence of sections
$\sigma_k$ which attains $\sup_{\mathcal{H}_1} J_{\alpha_k}$, and
analyze its blow-up behavior. First of all, we prove the following

\bthm{Lemma}{4.1.} {The functional $J_{\alpha_k}(\sigma)$ defined on
the space $\mathcal{H}_1$ admits a smooth maximizer
${\sigma}_k\in \mathcal{H}_1$. Moreover, we have
\begin{equation}\label{4.1}
\lim_{k\rightarrow+\infty}\int_Me^{\alpha_k|\sigma_k|^2}dV_g=
\sup_{{\sigma}\in\mathcal{H}_1}\int_Me^{4\pi |\sigma|^2}dV_g.
\end{equation}}\ethm

\proof It is easy to find  ${\sigma}_k\in
\mathcal{H}_1$ such that
$$J_{\alpha_k}({\sigma}_k)=\sup_{{\sigma}\in \mathcal{H}_1 }J_{\alpha_k}({\sigma}).$$
One can check that ${\sigma}_k$ satisfies the following Euler-Lagrange equation:
\begin{equation}\label{4.2}
    -\Delta {\sigma}_k=\frac{{\sigma}_k}{\lambda_k}
      e^{\alpha_k|\sigma_k|^2}-\sum_{i=1}^m\gamma^i_k\zeta_i,
\end{equation}
where $\Delta=\nabla^*\nabla$ is defined as follows: for any
$\tau\in H^{1,2}(M,E)$
$$\dint_M\lan\nabla\sigma_k,\nabla\tau\ran dV_g=\dint_M\lan\Delta \sigma_k,\tau\ran dV_g.$$
It is easy to check that
\begin{equation}\label{4.3}
      \lambda_k=\int_M |\sigma_k|^2
      e^{\alpha_k|\sigma_k|^2}dV_g,\s {\rm and} \s
    {\gamma}_k^i=\int_M \frac{\lan{\sigma}_k,\zeta_i\ran}{\lambda_k}
     e^{\alpha_k|\sigma_k|^2}dV_g.
\end{equation}

Let $\Omega_p$ be an coordinate system around $p$. We set
$$\sigma_k=u^i_ke_i,\x and\x U_k=(u_k^1,u_k^2\cdots,u_k^n),$$
and
$$\nabla e_j=\Gamma^l_{ji}dx^i\otimes e_l.$$
Then, we get a local version of (4.2)
\begin{equation}\label{4.4}
   -\Delta_0 u_k^j=T_{s}^{ji}\frac{\partial u^s_k}
    {\partial x^i}+R_{i}^ju_k^i
    +e^{f_p(x_k+r_kx)}(\frac{u_k^j}{\lambda_k}e^{\alpha_k|\sigma_k|^2}-\sum_{i=1}^m
    \gamma_k^i\lan\zeta_i,e_p\ran),
\end{equation}
where $T_{s}^{ji}$ and $R_{i}^j$ are smooth on $\Omega_p$.
(\ref{4.4}) can be written for simplicity
$$-\Delta_0U_k=
e^{f_p(x_k+r_kx)}\frac{U_k}{\lambda_k}e^{\alpha_k|U_k|^2}+T(\nabla
U_k)+R(U_k)- \sum_i\gamma_k^i\lan U_k,\zeta_i\ran.$$ The standard
elliptic estimates implies that ${\sigma}_k\in C^\infty(M)$.

Since for any fixed ${\sigma}\in\mathcal{H}$,
$$\int_Me^{\alpha_k|\sigma|^2}dV_g\leq \int_Me^{\alpha_k|\sigma_k|^2}dV_g,$$
(4.1) follows immediately.
$\hfill\Box$\\

We assume that
\begin{equation}\label{4.5}
   \lim_{k\rightarrow+\infty}\int_{M}e^{\alpha |\sigma_k|^2}dV_g=+\infty
\end{equation}
for all $\alpha>4\pi$ (Otherwise, by the weakly compactness of
$L^p$, passing to a subsequence, we have
$$\lim_{k\rightarrow+\infty}\int_{M}e^{\alpha |\sigma_k|^2}dV_g=
  \int_{M}e^{\alpha |\sigma_0|^2}dV_g,$$
where ${\sigma}_0$ is the weak limit of ${\sigma}_k$. Hence,
Theorem 1.1 holds). It follows from Lemma 3.5 that
$\lan\nabla\sigma_k,\nabla\sigma_k\ran dV_g
\rightharpoondown\delta_p$. Given any $\Omega\subset\subset
M\setminus\{p\}$, we take a cut-off function $\eta$ which is 0 at
$p$, and 1 on $\Omega$, then
$$\int_M\lan\nabla\eta\sigma_k,\nabla\eta\sigma_k\ran dV_g=\int_M\eta^2
\lan\nabla\sigma_k,\nabla\sigma_k\ran +2\eta\lan
\sigma_kd\eta,\nabla\sigma_k\ran dV_g+
\int_M\lan\sigma_kd\eta,\sigma_kd\eta\ran dV_g\rightarrow 0,$$
hence
$\int_\Omega(|\nabla\eta\sigma_k|^2+|\eta\sigma_k|^2)dV_g\rightarrow
0$. By Lemma 3.7, $e^{\alpha_k|\sigma_k|^2}$ is bounded in
$L^p(\Omega)$ for any $p>0$. Standard elliptic estimates imply
that $\sigma_k\rightarrow 0$ in $C^\infty(\Omega)$.

\bthm{Lemma}{4.2.} {Let $\frac{1}{\beta_k}=\int_M\frac{|\sigma_k|}{\lambda_k}
e^{\alpha_k|\sigma_k|^2}dV_g$. Then, we have
$$\liminf\limits_{k\rightarrow+\infty}\beta_k=+\infty,\,\,
\liminf\limits_{k\rightarrow+\infty}\frac{\lambda_k}{\beta_k}=+\infty,\,\,
 and\,\,\limsup\limits_{k\rightarrow+\infty}\gamma_k<+\infty.$$}\ethm

\proof For any fixed $N>0$, we have
$$\lim_{k\rightarrow+\infty}
\dint_{|\sigma_k|\leq N}e^{\alpha_k|\sigma_k|^2}dV_g=\dint_Me^0=Vol(M).$$
Then
$$\lim_{k\rightarrow+\infty}\dint_Me^{\alpha_k|\sigma_k|^2}\leq
Vol(M)+\lim_{k\rightarrow+\infty}\frac{1}{N^2}\dint_{|\sigma_k|>N}
|\sigma_k|^2e^{\alpha_k|\sigma_k|^2}dV_g
\leq Vol(M)+\lim_{k\rightarrow+\infty}\frac{\lambda_k}{N^2}.$$
If $\lim\limits_{k\rightarrow+\infty}\lambda_k<+\infty$, letting $N\rightarrow+\infty$,
we get
$$\int_Me^{\alpha_k|\sigma_k|^2}dV_g\leq Vol(M).$$
However, it follows from (\ref{4.1}) that
$$\lim_{k\rightarrow +\infty}\int_Me^{\alpha_k|\sigma_k|^2}dV_g
  =\sup_{{\sigma}\in\mathcal{H}_1}\int_Me^{4\pi |\sigma|^2}dV_g>Vol(M).$$
Therefore, $\liminf\limits_{k\rightarrow+\infty}
\lambda_k=+\infty.$

In the similar way, we have $\frac{\lambda_k}{\beta_k}\rightarrow
+\infty$. Then, we have
$$| \gamma^i_k|\leq c_1+c_2
   \int_{\{|\sigma_k|^2\geq 1\}}\frac{|\sigma|^2}{\lambda_k}
   e^{\alpha_k|\sigma_k|^2}dV_g\leq C
$$
for some constants $c_1$, $c_2$ and $C$ depending only on $M$.
Moreover, we have
$$\dint_M\frac{|\sigma_k|}{\lambda_k}e^{\alpha_k|\sigma_k|^2}dV_g
\leq \dint_M(\frac{e^{4\pi
N^2}}{\lambda_k}+\frac{1}{\lambda_kN}|\sigma|^2
e^{\alpha_k|\sigma_k|^2}dV_g).$$ Letting $k\rightarrow+\infty$,
and then $N\rightarrow+\infty$, we get
$\frac{1}{\beta_k}\rightarrow 0$.
  $\hfill\Box$\\

Let $c_k=|\sigma_k|(x_k)=\max_{x\in M}|\sigma_k|(x)$. By
(\ref{4.4}) and (\ref{4.5}), $c_k\rightarrow +\infty$. Passing to
a subsequence, we assume $x_k\rightarrow p$. Set $r_k^2=\lambda_k
c_k^{-2}e^{-\alpha_kc_k^2}$. Since
$$1=\int_M\frac{|\sigma_k|^2}{\lambda_k}
  e^{\alpha_k|\sigma_k|^2}dV_g\leq \frac{c_k^2}{\lambda_k}
  e^{\frac{\alpha_k}{2} c_k^2}\int_Me^{\frac{\alpha_k}{2}|\sigma_k|^2}dV_g
  \leq C\frac{c_k^2}{\lambda_k}
  e^{\frac{\alpha_k}{2} c_k^2},$$
we get
$$r_k^2 e^{\frac{\alpha_k}{2} c_k^2}\rightarrow 0,\,\,\, as\,\, k\rightarrow+\infty.$$

Let $\Omega_p$ be a local coordinate system around $p$,
$\sigma_k=u_k^ie_i$ and $U_k=(u_k^1,\cdots,u_k^n)$. Denote
$$V_k(x)=U_k(x_k+r_k x), \x D_k=U_k(x_k)$$
and
$$C_k=(c_k^1,c_k^2,\cdots,c_k^n)=
  \frac{1}{|B_1|}\int_{B_1}V_k dx.$$
A direct calculation shows
\begin{equation}\label{4.6}
  -\Delta_0 (V_k-C_k)=c_k^{-2}V_k
  e^{\alpha_k(|\sigma_k|^2-c_k^2)+f(x_k+r_kx)}+
  \mathcal{F}_k.
\end{equation}
where $||\mathcal{F}_k||_{L^2_{loc}}\rightarrow 0$. Hence
$$\begin{array}{lll}
  -\Delta_0(V_k-C_k)\rightarrow0\quad {\rm in}\quad
  L^2_{loc}(\mathbb{R}^2,\mathbb{R}^n), \x {\rm and}\x
  \dint_{B_L}|\nabla(V_k-C_k)|^2dx\leq 1 .
\end{array}$$
Then, by Poincare inequality, one has $\|V_k-C_k\|_{L^2(B_L)}\leq
C(L)$ for any fixed $L>0$. By the standard elliptic estimates,
$V_k-C_k\rightarrow V$ in
$C^{0,\alpha}_{loc}(\mathbb{R}^2,\mathbb{R}^n)$ for some $V\in
H^{1,2}(\mathbb{R}^2,\mathbb{R}^n)$. Notice that
$\int_{\mathbb{R}^2}\nabla_0V \nabla_0V^Tdx\leq 1$, Liouville's
theorem gives $V\equiv0$. Furthermore we obtain
\begin{equation}\label{4.7}
  V_k-D_k=(V_k-C_k)+(C_k-V_k(0))
  \rightarrow0 \,\,\,{\rm in}\,\,\, C^{0,\alpha}_{loc}(\mathbb{R}^2,\mathbb{R}^n).
\end{equation}
  Let
$$ w_k(x)=2\alpha_kD_k(U_k(x_k+r_k x)-D_k)^T.$$
Then we have
\begin{equation}\label{4.8}
  -\Delta_0  w_k=2\alpha_ke^{f(x_k+r_kx)}(1+\frac{D_k(V_k-D_k)^T}{c_k^2})
  e^{\alpha_k(|\sigma_k|^2(x_k+r_kx)-c_k^2)}+\mathcal{F}_k',
\end{equation}
where $||\mathcal{F'}_k||_{L^2_{loc}}\rightarrow 0$. It is easy to
see that on $B_L$
$$(D_kU_k^T)(x_k+r_kx)\leq |D_k|\times|U_k^T(x_k+r_kx)|
\leq c_k^2.$$
Hence
$$ w_k(x)\leq w_k(0)=0.$$
  Then the identity
$$\alpha_k(|\sigma_k|^2-c_k^2)=w_k+\alpha_k(V_k-D_k)(V_k-D_k)^T,$$
together with  Harnack inequality and standard elliptic estimates,
gives $ w_k\rightarrow w$ in $C^{0,\alpha}_{loc}(\mathbb{R}^2)$,
where $ w$ satisfies
$$\left\{\begin{array}{ll}
            -\Delta_0  w= 8\pi e^{ w}\,\,\,{\rm in}\,\,\, \mathbb{R}^2\\[1.2ex]
                w(0)=\sup_{\mathbb{R}^2} w=0\\[1.2ex]
            \int_{\mathbb{R}^2}e^{ w}dx\leq 1.
         \end{array}\right.$$
By a result of \cite{CL}, we have $ w(x)=-2\log(1+\pi |x|^2)$ in
$\mathbb{R}^2$ and $\int_{\mathbb{R}^2}e^{ w}=1$.\\

In the rest of this section, we will discuss the convergence of $\beta_k\sigma_k$.

\bthm{Proposition}{4.3.} $\beta_k {\sigma}_k\rightharpoondown G$ weakly in
$H^{1,q}(M,E)$ for any $1<q<2$, and $\beta_k {\sigma}_k\rightarrow G$ in
$C^2(\Omega)$ for any $\Omega\subset\subset M\setminus\{p\}$,
where $G\in C^2(M\setminus\{p\})$ satisfies
\begin{equation}\label{4.11}
   -\Delta G=\theta\delta_p-\sum_{i=1}^m\lan\theta,\zeta_i\ran\zeta_i
\end{equation}
for some $\theta\in E_p$.\ethm

\proof Since $\sigma_k\rightarrow 0$ in $C^\infty(M\setminus
B_\delta(p))$, and $\frac{\beta_k}{\lambda_k}\rightarrow 0$, then
for any $\varphi\in C^\infty(M,E)$, we have
$$\dint_{M\setminus B_\delta(p)}\frac{\beta_k\lan\varphi,\sigma_k\ran}{\lambda_k}e^{\alpha_k
|\sigma_k|^2}dV_g\rightarrow 0.$$
On $\Omega_p$, we set
$$\theta^i=\lim_{k\rightarrow+\infty}\dint_{\Omega_p}\frac{\beta_ku^i}{\lambda_k}e^{\alpha_k
|\sigma_k|^2}dV_g,$$
and $\theta=\theta^ie_i$. It is easy to see that
$$\dint_M\frac{\lan\varphi,\sigma_k\ran}{\lambda_k}e^{\alpha_k|\sigma_k|^2}dV_g\rightarrow
\lan\varphi,\theta\ran.$$
Moreover, we have
$$\dint_M\frac{\beta_k\lan\varphi,\sigma_k\ran}{\lambda_k}e^{\alpha_k|\sigma_k|^2}dV_g
=\dint_M\beta_k\gamma_k^i\lan\varphi,\zeta_i\ran dV_g$$
for any $\varphi\in C^\infty(M,E)$. Then, we get
$$|\dint_M\beta_k\gamma_k^i\lan\varphi,\zeta_i\ran dV_g|\leq ||\varphi||_{C^0},$$
therefore, $\sup\limits_k|\beta_k\gamma_k|<+\infty$.

Let $m_k^2=\int_M|\beta_k\sigma_k|^2dV_g$. Firstly, we need to
prove $\sup_km_k<+\infty$. If $m_k\rightarrow+\infty$, then
$||\Delta\frac{\beta_ku_k}{m_k}||_{L^1}\rightarrow 0$. It follows
from Proposition 7.1 in the Appendix that
$$\frac{\beta_k\sigma_k}{m_k}\rightharpoondown 0$$
in $H^{1,p}$ for any $p<2$. Hence it follows from Poincar\'{e}
inequality and compact embedding of Sobolev space that
$\frac{\beta_k\sigma_k}{m_k}\rightarrow 0$ in $L^2$, which
contradicts $\int_M |\frac{\beta_k\sigma_k}{m_k}|^2dV_g=1$.

Since $\sup_km_k<+\infty$, applying Proposition 7.1 again, we get
$\beta_k\sigma_k$ converges weakly in $H^{1,p}$ for any $p<2$.
Therefore $||\beta_k\sigma_k||_{L^q}\leq C(q)$ for any $q>0$.
Assume $\beta_k {\sigma}_k\rightharpoondown G$ weakly in
$H^{1,q}(M,E)$. Then we have by Lemma 4.3 and standard elliptic
estimate,
$$\beta_k {\sigma}_k\rightarrow G\quad\rm{in}\quad
C^2(\Omega),\quad\forall \Omega\subset\subset M\setminus\{p\}.$$
Testing the equation satisfied by $\beta_k\sigma_k$ with $\phi\in
C^\infty(M,E)$, one has
$$\begin{array}{lll}
     \int_M\lan\nabla\phi,\nabla \beta_k {\sigma}_k\ran dV_g&=&
      \int_M
     \frac{\lan\phi,\beta_k {\sigma}_k\ran}{\lambda_k}
     e^{\alpha_k|\sigma_k|^2}dV_g-\sum\limits_{i=1}^m\int_M\lan\phi,
     \beta_k\gamma_k^i\zeta_i\ran dV_g\\[1.7ex]
     &\rightarrow& \lan\theta,\phi(p)\ran-\sum_{i=1}^m\lan\theta,\zeta_i(p)\ran
     \int_M\lan\phi,\zeta_i\ran dV_g.
\end{array}$$
 Hence
$$\int_M\lan\nabla \phi, \nabla G\ran dV_g=\lan\theta,\phi\ran(p)-\sum_{i=1}^m
\lan\theta,\zeta_i(p) \ran\int_M\lan\phi,\zeta_i\ran dV_g,$$ and
whence (\ref{4.11}) holds. This completes the proof of the lemma.
$\hfill\Box$

\section{\bf Applications of Pohozaev identity}

In this section, we will calculate
$\lim\limits_{k\rightarrow+\infty}
\int_Me^{\alpha_k|\sigma_k|^2}dV_g$. The key gradient we use is
Pohozaev identity. Let $B_\delta(p)\subset\Omega_p$. Testing
equation (\ref{4.2}) with
$\eta(r)r\nabla_{\frac{\partial}{\partial r}} \sigma_k$, we have
$$\int_{B_\delta(p)}\lan\nabla \eta(r)r\nabla_{\frac{\partial}{\partial r}}\sigma_k,
\nabla\sigma_k\ran dV_g=\int_{B_\delta(p)} \lan\eta(r)r
\nabla_{\frac{\partial}{\partial r}}\sigma_k
,\sigma_k\ran\frac{e^{\alpha_k|\sigma_k|^2}}{\lambda_k}dx+\int_{B_\delta(p)}\gamma_k^i\lan
\eta(r)r\nabla_{\frac{\partial}{\partial r}}\sigma_k,\zeta_i \ran
dV_g,$$
where $\eta$ is a cut-off function which is 0 outside
$\partial B_\delta$ and 1 on $B_{\frac{\delta}{2}}$,
$r=\sqrt{(x^1)^2+(x^2)^2}$.
$$\begin{array}{lll}
    \lan\nabla \eta(r)r\nabla_{\frac{\partial}{\partial r}}\sigma_k,
      \nabla\sigma_k\ran_0&=&\sum\limits_{j}\lan\nabla_j\eta(r)x^p\nabla_p\sigma_k,
          \nabla_j\sigma_k\ran_0\\[1.7ex]
        &=&\eta(r)\lan\nabla\sigma_k,\nabla\sigma_k\ran_0
         +\eta'(r)\frac{x^ix^j}{r}\lan\nabla_i\sigma_k,\nabla_j\sigma_k\ran_0+
           \eta(r)\lan x^i\nabla_j\nabla_i\sigma_k,\nabla_j\sigma_k\ran_0\\[1.7ex]
        &=&\eta(r)\lan\nabla\sigma_k,\nabla\sigma_k\ran_0
           +\eta'(r)\frac{x^ix^j}{r}\lan\nabla_i\sigma_k,\nabla_j\sigma_k\ran_0+
          \eta(r)\lan x^i\nabla_i\nabla_j\sigma_k,\nabla_j\sigma_k\ran_0\\[1.7ex]
        &&+\eta(r)\sum\limits_{i,j}\lan x^i R_{ij}\sigma_k,\nabla_j\sigma_k\ran_0\\[1.7ex]
        &=&\eta(r)\lan\nabla\sigma_k,\nabla\sigma_k\ran_0
           +\eta'(r)\frac{x^ix^j}{r}\lan\nabla_i\sigma_k,\nabla_j\sigma_k\ran_0+
          \frac{1}{2}\eta(r)r\frac{\partial}{\partial r}\lan\nabla\sigma_k,
           \sigma_k\ran_0\\[1.7ex]
        &&+\eta(r)\sum\limits_{i,j}\lan x^i R_{ij}\sigma_k,\nabla_j\sigma_k\ran_0+\eta(r)
        S_0(\nabla\sigma_k,\nabla\sigma_k),
  \end{array}$$
where $S_0=e^{-f}S$, and
$$S(\sigma_1,\sigma_2)=\frac{1}{2}(r\frac{\partial}{\partial r}\lan \sigma_1,\sigma_2\ran
-\lan r\nabla_{\frac{\partial}{\partial r}}\sigma_1,\sigma_2\ran-
\lan\sigma_1, r\nabla_{\frac{\partial}{\partial
r}}\nabla\sigma_2\ran).$$
Since we do not
assume $\nabla h=0$, $S\neq 0$ generally.
We set
$$\eta(r)=\left\{\begin{array}{ll}
                      1&r\leq\delta\\
                      \frac{r-\delta-a}{-a}&\delta\leq r\leq \delta+a\\
                      0&r\geq\delta+a.
                 \end{array}\right.$$
Letting $a\rightarrow 0$, we get
$$\begin{array}{ll}
   \dint_{B_\delta(p)}\lan\nabla r\nabla_{\frac{\partial}{\partial r}}\sigma_k,
      \nabla\sigma_k\ran_0dx=&\dint_{\Omega_p}(\lan\nabla\sigma_k,\nabla\sigma_k\ran_0
           +\frac{1}{2}r\frac{\partial}{\partial r}\lan\nabla\sigma_k,\sigma_k\ran_0
           +\sum\limits_{i,j}\lan x^i R_{ij}\sigma_k,\nabla_j\sigma_k\ran_0)dx\\[1.7ex]
      &-\dint_{\partial B_\delta(p)}\frac{x^ix^j}{r}\lan\nabla_i\sigma_k,\nabla_j
      \sigma_k\ran_0dS+\dint_{B_\delta(p)}S(\nabla\sigma_k,\nabla\sigma_k)dV_g.
  \end{array}$$
Clearly,
$$\begin{array}{ll}
   \dint_{B_\delta(p)}\frac{1}{2}r\frac{\partial}{\partial r}
       \lan\nabla\sigma_k,\sigma_k\ran_0dx&=\dint_0^\delta \pi r^2\frac{\partial}{\partial r}
        (\dint_{S^1}\lan\nabla\sigma_k,\nabla\sigma_k\ran_0d\theta)dr\\[1.7ex]
     &=\frac{1}{2}\dint_{\partial B_\delta(p)}r\lan\nabla\sigma_k,\nabla\sigma_j\ran_0dS
       -\dint_{B_\delta(p)}\lan\sigma_k,\sigma_k\ran_0dx.
  \end{array}$$
Then
$$\begin{array}{ll}
     \dint_{B_\delta(p)}\lan\nabla r\nabla_{\frac{\partial}{\partial r}}\sigma_k,
       \nabla\sigma_k\ran_0dx=&\dint_{\Omega_p}
       \sum\limits_{i,j}\lan x^i R_{ij}\sigma_k,\nabla_j\sigma_k\ran_0dx
       -\dint_{\partial B_\delta(p)}\frac{x^ix^j}{r}\lan\nabla_i\sigma_k,\nabla_j
       \sigma_k\ran_0dS\\[1.7ex]
     &+\frac{1}{2}\dint_{\partial B_\delta(p)}r\lan\nabla\sigma_k,\nabla\sigma_j\ran_0dS
      +\dint_{B_\delta(p)}S(\nabla\sigma_k,\nabla\sigma_k)dV_g.
  \end{array}$$

On the other hand, we have
$$\dint_{B_\delta} e^{-f}\lan r\nabla_{\frac{\partial}{\partial r}}\sigma_k
       ,\sigma_k\ran_0\frac{e^{\alpha_k|\sigma_k|^2}}{\lambda_k}dx=\dint_{B_\delta}e^{-f}
    r\frac{\partial}{\partial r}\lan\sigma_k
       ,\sigma_k\ran_0\frac{e^{\alpha_k|\sigma_k|^2}}{\lambda_k}dx+\dint_{B_\delta}\frac{
       S(\sigma,\sigma)}{\lambda_k}
       e^{\alpha_k|\sigma_k|^2}dV_g$$
and
$$\begin{array}{lll}
    \dint_{B_\delta} e^{-f}\lan r\nabla_{\frac{\partial}{\partial r}}\sigma_k
       ,\sigma_k\ran_0\frac{e^{\alpha_k|\sigma_k|^2}}{\lambda_k}dx
       &=&\dint_0^\delta2\pi r^2\frac{\partial}{\partial r}(\dint_{S^1}
       \frac{e^{\alpha_k|\sigma|^2}}{2\alpha_k\lambda_k}e^{-f}d\theta)dr
       +\dint_{B_\delta}\frac{e^{\alpha_k|\sigma|^2}}{\alpha_k\lambda_k}r
        \frac{\partial f}{\partial r}dV_g.\\[1.7ex]
       &=&\dint_{\partial B_\delta(p)}r\frac{e^{\alpha_k|\sigma_k|^2}}
           {\alpha_k\lambda_k}dS_g-\dint_{B_\delta(p)}\frac{e^{\alpha_k|\sigma_k|^2}}{\alpha_k
           \lambda_k}dV_g\\[1.7ex]
       &&+\dint_{B_\delta}\frac{e^{\alpha_k|\sigma|^2}}{\alpha_k\lambda_k}r
        \frac{\partial f}{\partial r}dV_g.
  \end{array}$$
Therefore, we get
$$\begin{array}{lll}
   \dint_{B_\delta(p)}(1+O(r))\frac{e^{\alpha_k|\sigma_k|^2}}{\alpha_k\lambda_k}dV_g
     &=&\dint_{\partial B_\delta(p)}r\frac{e^{\alpha_k|\sigma_k|^2}}
           {\alpha_k\lambda_k}dS_g-\dint_{B_\delta(p)}
       \sum\limits_{i,j}\lan x^i R_{ij}\sigma_k,\nabla_j\sigma_k\ran_0dx\\[1.7ex]
     &&-\dint_{\partial B_\delta(p)}\frac{x^ix^j}{r}\lan\nabla_i\sigma_k,\nabla_j
       \sigma_k\ran_0dS-\frac{1}{2}\dint_{\partial B_\delta(p)}
       r\lan\nabla\sigma_k,\nabla\sigma_k\ran_0dS\\[1.7ex]
     &&+\dint_{B_\delta}\gamma_k^i\lan
        \eta(r)r\nabla_{\frac{\partial}{\partial r}}\sigma_k,\zeta_i \ran dV_g\\[1.7ex]
     &&+\dint_{B_\delta}(S(\nabla\sigma_k,\nabla\sigma_k)-\frac{S(\sigma_k,\sigma_k)}
      {\lambda_k}e^{\alpha_k|
     \sigma_k|^2})dV_g.
  \end{array}$$
Firstly, letting $k\rightarrow+\infty$, we have
$$\dint_{B_\delta}\gamma_k^i\lan
        \eta(r)r\nabla_{\frac{\partial}{\partial r}}\beta_k\sigma_k,\zeta_i \ran dV_g
\leq C||\nabla \beta_k\sigma_k||_{H^{1,q}_{B_\delta}}\rightarrow
C||\nabla G||_{H^{1,q}(B_\delta)}$$
and
$$\dint_{B_\delta(p)}
       \sum\limits_{i,j}(\lan x^i R_{ij}\beta_k\sigma_k,\nabla_j\beta_k\sigma_k\ran_0+S)dx
   \leq C||\beta_k\sigma_k||_{L^{q'}(B_\delta)}||\nabla\beta_k\sigma_k||_{L^q(B_\delta)
   }\rightarrow
   C||G||_{L^{q'(B_\delta)}}||\nabla G||_{L^q(B_\delta)},$$
where $1<q<2$. Secondly, applying Lemma 7.2 in the Appendix, we get
$$\begin{array}{l}
   \dint_{\partial B_\delta(p)}\frac{x^ix^j}{r}\lan\nabla_i\beta_k\sigma_k,\nabla_j
       \beta_k\sigma_k\ran_0dS_0+\frac{1}{2}\dint_{\partial B_\delta(p)}
       r\lan\nabla\beta_k\sigma_k,\nabla\beta_k\sigma_j\ran_0dS\\[1.7ex]
    \rightarrow
       \dint_{\partial B_\delta(p)}\frac{x^ix^j}{r}\lan\nabla_iG,\nabla_j
       G\ran_0dS-\frac{1}{2}\dint_{\partial B_\delta(p)}
       r\lan\nabla G,\nabla G\ran_0dS\\[1.7ex]
    =\frac{1}{4\pi}+O(\delta^{\gamma}).
  \end{array}$$
Thirdly, since $\sigma_k\rightarrow 0$ in $\partial B_\delta$, we get
$$\dint_{\partial B_\delta(p)}\frac{e^{\alpha_k|\sigma_k|^2}}{\alpha_k}
\rightarrow \frac{|\partial B_\delta|}{4\pi}.$$ Finally, we set
$S(\sigma_1,\sigma_2)=\lan A\sigma_1,\sigma_2\ran$, where
$||A||_{C^0}=O(r)$. Take a cut-off function $\eta$ which is 0
outsider $B_{2\delta}$ and 1 in $B_\delta$ with $|\nabla\eta|<\frac{1}{\delta}$. We have
$$\begin{array}{ll}
    \dint_{B_{2\delta}(p)}\lan A\nabla\eta\sigma_k,\nabla\sigma_k\ran
     &=\dint_{B_{2\delta}(p)}\lan\nabla A\eta\sigma_k,\nabla\sigma_k\ran dV_g
      -\dint_{B_{2\delta}(p)}\lan(\nabla A)\sigma_k,\nabla\sigma_k\ran dV_g\\[1.7ex]
     &=\dint_{B_{2\delta}(p)}\frac{\lan A\eta\sigma_k,\sigma_k\ran} {\lambda_k} e^{\alpha_k
     |\sigma_k|^2}dV_g
       -\dint_{B_{2\delta}(p)}\lan(\nabla A)\sigma_k,\nabla\sigma_k\ran dV_g
  \end{array}$$
Hence, we get
$$\begin{array}{l}
   \beta_k^2\dint_{B_\delta}(S(\nabla\sigma_k,\nabla\sigma_k)-S(\sigma_k,\sigma_k)\frac{e^{\alpha_k
     |\sigma_k|^2}}{\lambda_k})dV_g\\
   \x\x\x\x\rightarrow O(1)\dint_{B_{2\delta}(p)}\lan(\nabla A)\eta G,\nabla G\ran dV_g+
     O(\delta)\dint_{B_{2\delta}(p)\setminus B_\delta(p)}\eta|\nabla G|^2dV_g.
\end{array}$$
So, we get
$$\lim_{k\rightarrow+\infty}\dint_{B_\delta(p)}(1+O(r))e^{\alpha_k|\sigma_k|^2}dV_g=
\frac{|\partial B_\delta|}{4\pi}+\lim_{k\rightarrow+\infty}\frac{\lambda_k}{\beta_k^2}
+h(\delta),$$
where $\lim_{\delta\rightarrow 0}h(\delta)=0$.
Clearly, $\lim\limits_{k\rightarrow+\infty}\int_{M\setminus B_\delta(p)}
e^{\alpha_k|\sigma_k|^2}dV_g
=|M\setminus B_\delta(p)|$, we get
\begin{equation}{\label{5.1}}\lim_{k\rightarrow+\infty}\dint_Me^{\alpha_k|\sigma_k|^2}dV_g
=Vol(M)+\lim_{k\rightarrow+\infty}
\frac{\lambda_k}{\beta_k^2}.
\end{equation}

The following lemma is very important:

\bthm{Lemma}{5.1.} {we have $\frac{c_k}{\beta_k}\rightarrow 1$}.\ethm

\proof Since
$$\frac{c_k}{\beta_k}=\int_M\frac{c_k\sigma_k}
 {\lambda_k}e^{\alpha_k|\sigma_k|^2}\geq \int_M\frac{|\sigma_k|^2}{\lambda_k}
 e^{\alpha_k|\sigma_k|^2}dV_g,$$
we have $c_k\geq \beta_k$.

 Assume the lemma is not true, then we can find $A>1$ such that
 $\frac{c_k}{A}>A\beta_k$ for sufficiently large $k$. Denote
$$\tau=\lim_{k\rightarrow+\infty}\frac{\sigma_k}{c_k}=\tau^ie_i.$$
Without loss of generality, we assume $\tau^1> 0$. By the equation
(\ref{4.4}), we have  for any fixed $L>0$,
$$\begin{array}{lll}
 \dint_{M}|\nabla_0\eta(u_k^i-
 \frac{u_k^i(x_k)}{A})^+|^2dx&=&\dint_{\Omega_p}\nabla_0 u_k^i\nabla_0\eta(u_k^i-
    \frac{u_k^i(x_k)}{A})^+dx+o_k(1)\\[1.7ex]
 &=&\dint_{\Omega_p}\eta(u_k^i-
    \frac{u_k^i(x_k)}{A})^+\frac{u_k^i}{\lambda_k}
    e^{\alpha_k|\sigma_k|^2}dV_g+o_k(1)\\[1.7ex]
 &\geq&\dint_{B_{Lr_k(x_k)}}
    \frac{(u_k^i-
    \frac{u_k^i(x_k)}{A})u_k^i}{\lambda_k}
   e^{\alpha_k|\sigma_k|^2}dV_g+o_k(1)\\[1.7ex]
 &=&\dint_{B_L(0)}\left(\frac{u_k^i(x_k)(A-1)+o_k(1)}{A}\right)
    \frac{u_k^i}{c_k^2}
   e^{\varphi_k(x)+o_k(1)}dx+o_k(1).
 \end{array}$$
Letting $k\rightarrow+\infty$, then $L\rightarrow+\infty$, we get
 $$\liminf_{k\rightarrow+\infty}\int_M|\nabla\eta (u_k^1-
 \frac{u_k^1(x_k)}{A})^+|^2dV_g
 \geq (1-\frac{1}{A}){\tau^1}^2.$$
Let $\sigma_k^A=\eta
u_k^{A,i}e_i=\eta(\min\{u_k^1(x),\frac{u_k(x_k)}{A}\}e_1+
\sum\limits_{i=2}^n u_k^ie_i)$, where $\eta$ is a cut-off function
which is 0 $M\setminus B_{2\delta}(p)$, and 1 in $B_\delta(p)$.
Since
$$|\nabla\sigma_k^A|^2=\eta\sum_{i=1}|\nabla_0u_k^{A,i}|^2+O(|\sigma_k^A||\nabla\sigma_k^A|)
+O(|\sigma_k^A|^2),$$
we get
$$\limsup_{k\rightarrow+\infty}\dint_M(|\nabla \sigma_k^A|^2+|\sigma_k^A|^2)dV_g\leq
1-\frac{A-1}{A}\tau^{1^2}<1.$$
Then, by Proposition 3.7,
$e^{\alpha_k|\sigma_k^A|}$ is bounded in $L^q$ for some $q>1$.
Since
$$\dint_{\Omega_p\cap \{u_k^1\leq \frac{\tau^1c_k}{A}\}}
\frac{\beta_k}{\lambda_k}|\sigma_k|e^{\alpha_k|\sigma_k|^2}dV_g
\leq
\frac{\beta_k}{\lambda_k}\dint_{\Omega_p}|\sigma_k|e^{\alpha_k|\sigma_k^A|^2}dV_g
\leq
\frac{\beta_k}{\lambda_k}||\sigma_k||_{L^{q'}(M)}||e^{\alpha_k|\sigma_k^A|^2}||_{L^q(M)}
\rightarrow 0.$$
In the same way, for any $\tau^i$, we have
$$\dint_{\Omega_p\cap \{|u_k^i|\leq \frac{\tau^ic_k}{A}\}}
\frac{\beta_k}{\lambda_k}|\sigma_k|e^{\alpha_k|\sigma_k|^2}dV_g\rightarrow 0$$

Since $|\tau|=1$, we obtain
$$\begin{array}{ll}
   1&\geq\lim\limits_{k\rightarrow+\infty}\dint_{\{|\sigma_k|\geq A\beta_k\}}\frac{
   |\sigma_k|^2}
       {\lambda_k}e^{\alpha_k|\sigma_k|^2}dV_g\\[1.7ex]
     &\geq A\lim\limits_{k\rightarrow+\infty}\beta_k\dint_{\{x\in\Omega_p:|u_k^i|\geq A\tau^i
     \beta_k,\forall i \}}
        \frac{|\sigma_k|}
       {\lambda_k}e^{\alpha_k|\sigma_k|^2}dV_g\\[1.7ex]
     &\geq A\lim\limits_{k\rightarrow+\infty}\beta_k\dint_M\frac{|\sigma_k|}
       {\lambda_k}e^{\alpha_k|\sigma_k|^2}dV_g-A\lim_{k\rightarrow+\infty}\sum\limits_{i}
      \dint_{\{|u_k^i|\leq \frac{\tau^ic_k}{A}\}}\frac{\beta_k|\sigma_k|}{\lambda_k}
      e^{\alpha_k|\sigma_k|^2}dV_g\\[1.7ex]
     &=A.
  \end{array}$$
This contradict with the choice of $A$.
$\hfill\Box$\\

\bthm{Corollary}{5.2.} We have $\tau=\theta$, and
\begin{equation}\label{5.2}
    \lim_{k\rightarrow+\infty}\dint_M(e^{\alpha_k|\sigma_k|^2}-1)dV_g=
    \lim_{L\rightarrow+\infty}\lim_{k\rightarrow+\infty}\dint_{B_{Lr_k}(x_k)}
    e^{\alpha_k|\sigma_k|^2}dV_g.
\end{equation}\ethm

\proof  By a straightforward calculation, we have
$$\dint_{B_{Lr_k}(x_k)}e^{\alpha_k|\sigma_k|^2}
dV_g=\frac{\lambda_k}{c_k^2}\dint_{B_L}e^wdx(1+o(1)).$$ The above
inequality together with (\ref{5.1}) implies (\ref{5.2}).

It is not difficult to check that
$$\lim_{L\rightarrow+\infty}\lim_{k\rightarrow+\infty}\dint_{B_{Lr_k}(x_k)}
\frac{\beta_k|\sigma_k|}{\lambda_k}e^{\alpha_k|\sigma_k|^2}dV_g=1,$$
and
$$\lim_{L\rightarrow+\infty}\lim_{k\rightarrow+\infty}\dint_{B_{Lr_k}(x_k)}
\frac{\beta_ku_k^i}{\lambda_k}e^{\alpha_k|\sigma_k|^2}dV_g=\tau^i.$$
Hence
$$|\theta^i-\tau^i|=\lim_{L\rightarrow+\infty}\lim_{k\rightarrow+\infty}
\int_{\Omega_p\setminus B_{Lr_k}(x_k)}
\frac{\beta_ku_k^i}{\lambda_k}e^{\alpha_k|\sigma_k|^2}dV_g \leq
1-\lim_{L\rightarrow+\infty}\lim_{k\rightarrow+\infty}
\int_{B_{Lr_k}(x_k)}
\frac{\beta_k|\sigma_k|}{\lambda_k}e^{\alpha_k|\sigma_k|^2}dV_g=0.$$
$\hfill\Box$\\

\section{\bf The proof of theorem 1.1}

On $\Omega_p$, we set
$\tilde{U}_k=(|\tilde{u}_k^1|,\cdots,|\tilde{u}_k^n|)$. The
following lemma is very important for the rest of our arguments:

\bthm{Lemma}{6.1.} We have
$$\int_{B_\delta(p)}|\nabla_0 \tilde{U}_k|^2dx\leq\int_{B_\delta(p)}
\lan\nabla {\sigma}_k,\nabla\sigma_k\ran dV_g+\frac{\rho(\delta)}{c_k^2},$$ where
$\rho(\delta)$ is a continuous function of $\delta$ with $\rho(0)=0$.\ethm

\proof It is well-known that
$$\int_{B_\delta(p)}|\nabla_0\tilde{U}_k|^2dx=\int_{B_\delta(p)}|\nabla_0 U_k|^2dx.$$
We set
$$\nabla {\sigma}_k=(\frac{\partial u_k^p}{\partial x^1}+
\frac{\partial u_k^p}{\partial x^2})e_p+A^p_iu_k^ie_p,$$ where
$A^p$'s are smooth functions. Hence, we get
$$|\nabla_0U_k|^2=|\nabla_0U_k|^2-2\lan A(\sigma_k),\nabla\sigma_k\ran+|AU_k|^2.$$
Since
$$\begin{array}{ll}
 \dint_{B_\delta(p)}c_k^2\lan A(\sigma_k),\nabla\sigma_k\ran dx
   &\leq C(\int_{B_\delta(p)}|A(c_k {\sigma}_k)|^qdx)^\frac{1}{q}(\int_{B_\delta(p)}
      |\nabla c_k {\sigma}_k|^pdx)^\frac{1}{p}\\[1.7ex]
   &\leq C(\int_{B_\delta(p)}|c_k {\sigma}_k|^qdx)^\frac{1}
        {q}(\int_{B_\delta(p)}
        |\nabla c_k {\sigma}_k|^pdx)^\frac{1}{p}\\[1.7ex]
   &\rightarrow C(\int_{B_\delta(p)}|G|^qdx)^\frac{1}
        {q}(\int_{B_\delta(p)}
        |\nabla G|^pdx)^\frac{1}{p},
\end{array}$$
we get
$$\int_{B_\delta(p)}\lan A({\sigma}_k), \nabla \sigma_k\ran dV_g\leq \frac{\rho_1(\delta)}
{c_k^2}.$$
Similarly we have
$$\int_{B_\delta(p)}|A(\sigma_k)|^2\leq \frac{\rho_2(\delta)}{c_k^2}.$$
The lemma follows immediately from the above two inequalities.
$\hfill\Box$\\

On $\Omega_p$, we can write $G$ as follows :
$$G=-\frac{\log r}{2\pi}{\tau}+s_p+\theta.$$
where $s_p$ is a constant section, and $\theta$ is continuous
local section of $E$ with $\theta(0)=0$. Then Theorem 1.1 follows
from the following proposition:

\bthm{Proposition}{6.2} If $(\ref{4.5})$ holds, then we have
$$\sup_{{\sigma}\in \mathcal{H}_1}\int_Me^{4\pi |\sigma|^2}dV_g\leq Vol(M)+\pi
 e^{1+4\pi\lan\tau,s_p\ran}.$$\ethm

\proof Set $Y_k=(\tilde{U}_k-L_k)^+$, where
$$L_k=(l_1,\cdots,l_n)=\sup_{\partial B_\delta(p)}\tilde{U}_k=O(\frac{1}{c_k}).$$
Then $Y_k\in H^{1,2}_0(B_\delta,\mathbb{R}^n)$, and
$$\dint_{B_\delta(p)}|\nabla_0 Y_k|^2dx\leq 1-\int_{M\setminus B_\delta(p)}|\nabla G|^2dV_g
+o_\delta(1)\frac{1}{c_k^2}.$$
By Lemma 6.1 and  Proposition 7.2
in Appendix, we have
$$\int_{B_{\delta}(p)}|\nabla Y_k|^2dx\leq
1-\frac{-\frac{1}{2\pi}\log{\delta}+ \lan\tau,s_p\ran
+\epsilon(k,\delta)}{c_k^2},$$ where
$\lim\limits_{\delta\rightarrow+0}\lim\limits_{k\rightarrow+\infty}
\epsilon(k,\delta)=0$. Let
$$\vartheta_k=\frac{1}{1-\frac{-\frac{1}{2\pi}\log{\delta}+
\lan\tau,s_p\ran +\epsilon(k,\delta)}{c_k^2}}=
1+\frac{-\frac{1}{2\pi}\log{\delta}+ \lan\tau,s_p\ran
+\epsilon(k,\delta)}{c_k^2}+O(\frac{\log^2{\delta}}{c_k^4}).$$
By
Theorem A and the method we used to prove Corollary 3.2, we have
$$\limsup_{k\rightarrow+\infty}\dint_{B_\delta(p)}e^{\alpha_k\vartheta_k|Y_k|^2}dx
\leq \delta^2(\pi e+\pi).$$

If $\tilde{u}_i<l_k$, we have on $B_{Lr_k}(x_k)$
$$\tau^i=0,\quad \tilde{u}_i=\frac{o_\delta(1)}{c_k},\quad l_k=\frac{o_\delta(1)}{c_k}.$$
Hence, we have $y_k^i=\tilde{u}_k^i-l_k+\frac{o_\delta(1)}{c_k}$
for sufficiently large $k$. A straightforward calculation shows on
$B_{Lr_k}(x_k)$,
$$\vartheta_k|Y_k|^2=|\tilde{U}_k|^2
       -2\tilde{U}_kL_k+|\tilde{U}_k|^2(
       \frac{-\frac{1}{2\pi}\log{\delta}+
\lan\tau,s_p\ran}{c_k^2})+\epsilon(\delta,k).$$
It is easy to check that
$$\lim_{k\rightarrow+\infty}c_kL^i_k=-\frac{|{\tau^i}|}{2\pi}\log\delta+sign({\tau^i})s_p^i
+o_\delta(1),$$ and on $B_{Lr_k}(x_k)$,
$$\frac{\tilde{U}_k}{c_k}=(|\tau^1|,\cdots,|\tau^n|)+\frac{o_k(1)}{c_k}.$$
Recall that $|\tau|=1$, we get
$$\tilde{U}_kL_k=-\frac{1}{2\pi}\log\delta+\lan\tau,s_p\ran+
\epsilon'(\delta,k),$$
where $\lim\limits_{\delta\rightarrow 0}\lim\limits_{k\rightarrow+\infty}
\epsilon'(\delta,k)=0$.
Since $|\tilde{U}_k|^2/c_k^2\rightarrow1$ as $k\rightarrow+\infty$, we have
$$|\tilde{U}_k|^2\leq \vartheta_k|Y_k|^2
-\frac{1}{2\pi}\log{\delta}+
\lan\tau,s_p\ran+\epsilon'(\delta,k).$$ Letting
$k\rightarrow+\infty$, then $L\rightarrow+\infty$, and
$\delta\rightarrow 0$, we get
\begin{equation}\label{6.1}
  \lim_{L\rightarrow+\infty}\lim_{k\rightarrow+\infty}\int_{B_{Lr_k(x_k)}}
  e^{\alpha_k|\sigma_k|^2} dV_g \leq \pi e^{1+4\pi\lan\tau,s_p\ran}
  \lim_{L\rightarrow+\infty}\lim_{k\rightarrow+\infty}\int_{B_{Lr_k(x_k)}}
   e^{\alpha_k\vartheta_k|Y_k|^2} dx<+\infty.
\end{equation}
Then
by Corollary 5.2, we obtain
$$\lim_{L\rightarrow+\infty}
\lim_{k\rightarrow+\infty}\dint_{B_\rho(p)\setminus B_{Lr_k}(x_k)}
e^{\alpha_k|\sigma_k|^2}dV_g=|B_\rho|,\quad \forall\rho>0.$$
Clearly,
$$|Y_k|^2\leq |\tilde{U}_k|^2=|\sigma_k|^2,\x and\x
\alpha_k\vartheta_k|Y_k|^2\leq \alpha_k|U_k|^2-2\log{\delta}+C.$$
Hence
$$\lim_{L\rightarrow+\infty}\lim_{k\rightarrow+\infty}
\dint_{B_\rho(p)\setminus B_{Lr_k}(x_k)}
 e^{\alpha_k\vartheta_k|Y_k|^2}dx\leq
 O(\delta^{-2})\dint_{B_\rho(p)}e^{\alpha_k|\sigma_k|^2}dV_g=O(\delta^{-2})|B_\rho(p)|.$$
It is easy to see that for any fixed $\rho>0$
$$\dint_{B_\delta(p)\setminus B_\rho(p)}e^{\alpha_k\vartheta_k|Y_k|^2}dV_g\rightarrow
|B_\delta\setminus B_\rho|.$$
Letting $\rho\rightarrow 0$, we get
$$\lim_{L\rightarrow+\infty}
  \lim_{k\rightarrow+\infty}\dint_{B_\delta(p)\setminus B_{Lr_k}(x_k)}e^{\alpha_k\vartheta_k|Y_k|^2}dx
  =\pi\delta^2,\lim_{L\rightarrow+\infty}
  \lim_{k\rightarrow+\infty}\dint_{B_{Lr_k}(x_k)}e^{\alpha_k\vartheta_k|Y_k|^2}dx
  \leq \delta^2\pi e.$$
Then the Proposition
follows from (\ref{6.1}) and Corollary 5.2. $\hfill\Box$

\section{\bf Appendix}
{\small We will prove two propositions which have been used in
section 5 and section 6. Since the proof is routine, we prove them
in this appendix. The first proposition is an extension of Theorem
2.2 in \cite{S2}:

\bthm{Proposition}{7.1.} Let $\sigma$ be a  smooth section of $E$ with the equation
$$\Delta\sigma=f.$$
If $||\sigma||_{L^2}\leq \gamma$ and $||f||_{L^1}\leq 1$, then for any $q<2$,
there is a constant $C(q)$ which
depend only on $q$ $\gamma$ and $M$, s.t.
$$||\sigma||_{H^{1,q}}\leq C(q)$$\ethm

{\it Proof.} For any $p\in M$, we take a cut-off function which is
1 in $B_r$ and $0$ outside $B_{2r}$. Write $\sigma=u^ie_i$. Given
$t>0$, we set $u^{i,t}=\min\{\eta u^i,t\}$ and
$\sigma^t=u^{i,t}e_i$. Then
$$\dint_M\lan f,\eta^2\sigma^t\ran dV_g=\int_M\lan\eta^2\sigma^t,\Delta\sigma\ran
dV_g.$$
We get
$$\dint_{B_r(p)}\sum_{i}|\nabla\eta u^{i,t}|^2dx\leq C_1t+C_2\leq C't$$
when $t$ is sufficiently large.

Let $u^i_*$ be the  rearrangement of $\eta u^i$, and
$\mu_{\mathbb{R}^2}(B_\rho)=\mu_{\mathbb{R}^2}\{ \eta u^{i,t}\geq
t\}$. Then we have
$$\inf\left\{\dint_{B_{r}}|\nabla v|^2dx: v\in H^{1,2}_0(B_{r})
 \hbox{ and } v|_{B_\rho}=t\right\}\leq Ct.$$
On the other hand, the inf is achieved by
$-t\log{\frac{|x|}{r}}/\log{ \frac{r}{\rho}}$. By a direct
computation, we have
$$\frac{2\pi t}{\log\frac{r}{\rho}}\leq C ,$$
$$\mu_{M}(\{x\in B_{r}:\eta u^i\geq t \})\leq C\mu_{\mathbb{R}^n}(\{x\in B_{r}:\eta u^i\geq t
 \})
 =C\mu_{\mathbb{R}^n}(B_\rho)\leq C(r,p)e^{-A(r,p)t}.$$
Hence, we can find a constant $A$, s.t.
$$|\{\eta u^i\geq t\}\cap B_r|\leq e^{-A t}.$$
In the same way, we get
$$|\{\eta u^i\leq -t\}\cap B_r|\leq e^{-A' t}.$$
By the compactness of
$M$, we get
$$|\{|\sigma|\geq t\}|\leq e^{-\delta_0 t}$$
for some $\delta_0$ and sufficiently large $t$. Then, for any
$\delta<\delta_0$,
$$\dint_Me^{\delta |\sigma|}dV_g\leq \sum_{m=0}^\infty\mu(\{m\leq u
\leq m+1\})e^{\delta(m+1)}\leq \sum_{m=0}^\infty
e^{-(\delta_0-\delta)m}e^\delta\leq C.$$ Locally, we denote
$U=(u^1,\cdots,u^n)$, then
$$\Delta_0u^i=f^i+T^i_1(\nabla U)+T_2^i(U),$$
where $T_1$ and $T_2$ are smooth and linear at each point. Testing
the above equation with the function
$\eta^2\log{\frac{1+2{u^i}^+}{1+{u^i}^+}}$, where $\eta$ is 1 on
$B_r$ and $0$ outside $B_{2r}$. We obtain the following
$$\int_{B_r}\sum_{i=1}^n\frac{|\nabla_0\eta {u^i}^+|^2}{(1+2{u^i}^+)(1+{u^i}^+)}
dx\leq
\log{2}+\int_{B_r}(T_{1'}^i(\nabla\eta{U}^+)+T_{2'}^i({U}^+))dV_g
,$$ where $T_{1'}$ and $T_{2'}$ are smooth and linear at each
point. Since $\int_M e^{\delta|\sigma|}dV_g<C(\delta)$, we have
$$\int_{B_r}\sum_{i=1}^n\frac{|\nabla {u^i}^+|^2}{(1+2{u^i}^+)(1+{u^i}^+)}\leq C.$$
Given $q<2$, by Young's inequality, we have
$$\begin{array}{ll}
   \dint_{B_r}|\nabla {u^i}^+|^qdV_g&\leq\dint_{B_r}\left(\frac{|\nabla_0 {u^i}^+|^2}
   {(1+{u^i}^+)(1+2{u^i}^+)}+((1+{u^i}^+)(
            1+2{u^i}^+))^\frac{q}{2-q}\right)dx\\[1.7ex]
            &\leq \dint_{B_r}\left(\frac{|\nabla_0
              {u^i}^+|^2}{(1+{u^i}^+)(1+2{u^i}^+)}+Ce^{\delta {u^i}^+}\right)dx.
            \end{array}$$
In the same way,
$$\dint_{B_r}|\nabla {u^i}^-|^qdV_g<C(q).$$
Again, by the compactness of $M$, we get
$$||\sigma||_{H^{1,q}}<C(q)$$
for some $C(q)$.

$\hfill\Box$\\

On $\Omega_p$, we can set
$$G=-\frac{\log r}{2\pi}{\tau}+s_p+\theta,$$
we have the following Lemma:\\

\bthm{Lemma}{7.2.} There are constants $\gamma\in (0,1)$ and
$A>0$, s.t.
$$|\nabla \theta|(x)\leq Ar^{\gamma-1},$$
when $x$ is near 0. Therefore, we have
$$\int_{M\setminus B_\delta(p)}|\nabla G|^2dV_g=-\frac{1}{2\pi}\log{\delta}
+\lan\tau,s_p\ran+o_\delta(1).$$\ethm

\proof Let $\theta=v^ie_i$, and $V=(v^1,\cdots,v^n)$,
we have the equation of $V$
$$-\Delta_0 V=-\nabla Q_1\nabla V+\frac{1}{r}Q_2\frac{\partial Q_3}
{\partial r}+F,$$
where $Q_1$ $Q_2$ and $Q_3$ are smooth matrix, $F$ is a smooth vector function
. Hence $V\in W^{2,p}_{loc}$ for any $p<2$.

Clearly, we can write
$$Q^{-1}\frac{\partial Q}{\partial r}=
Q^{-1}\frac{\partial Q}{\partial x_i}\frac{x_i}{r}=
\Lambda_i\frac{x_i}{r}+F_1,$$
where $\Lambda_k$ are constant vector, and $F_1=o(1)$. Hence
$$-\Delta_0 V=\Lambda_i\frac{x_i}{r^2}+F_2,$$
where $F_2\in L^\infty_{loc}$. Then
$$-\Delta_0 x_mV=\Lambda_i\frac{x_ix_m}{r^2}+x_mF_2+\frac{\partial V}{\partial x_m}
\in L^q_{loc}$$
for any $q>0$. Hence, we get
$$||x_mV||_{C^1(B_r)}\leq ||V||_{C^0(B_{2r})}2r+C(1+||\nabla V||_{L^{2q}})
r^{\frac{1}{q'}}$$
for any $q>1$.
Notice that $V(0)=0$, and then $V=O(r^\gamma)$ for some $\gamma>0$, we get
$$|x_m\frac{\partial V}{\partial x_i}(x)|\leq C r^{\gamma}$$
for some $\gamma>0$ , any $i,m=1,2$ and $x\neq 0$. In the end, we get
$$|\nabla V|(x)\leq |x|^{\gamma-1}.$$
     $\hfill\Box$\\
}

\bigskip

\bigskip
\bigskip

\noindent Yuxiang  Li\\
{ICTP, Mathematics Section, Strada Costiera 11, I-34014 Trieste,
Italy}\\
{\it E-mail address:} liy@ictp.it\\

\noindent Pan Liu\\
{Department of Mathematics, East China Normal University, 3663, Zhong Shan North Rd,
Shanghai 200062, P.R.China\\
{\it E-mail address:} pliu@math.ecnu.edu.cn\\

\noindent Yunyan Yang\\
{Department of Mathematics, Information School, Renmin University
of China, Beijing 100872, P. R. China}\\ {\it E-mail address:}
yunyan$\_$yang2002@yahoo.com.cn

\end{document}